\documentclass[12pt]{article}
\usepackage[dvipsnames]{xcolor}
\usepackage{amsmath}
\usepackage{setspace}
\usepackage{bm}
\usepackage{amssymb}
\usepackage{indentfirst}
\usepackage[superscript]{cite}
\usepackage{geometry}
\usepackage{bbm}
\usepackage{tocloft}
\usepackage{mathtools}
\definecolor{mygreen}{RGB}{0,128,0}

\DeclareMathOperator{\HP}{HP}

\newcommand{\ii}{\bm{i}}

\usepackage[super]{natbib}
\setcitestyle{super}

\newcommand*{\citena}[1]{%
\begingroup
[\color{Green}
\romannumeral-`\x % remove space at the    beginning of \setcitestyle
\setcitestyle{numbers}%
\cite{#1}%
\endgroup
]\ignorespacesafterend
}

\newcommand*{\citesup}[1]{%
\begingroup
\color{Green}
\cite{#1}%
\endgroup
\ignorespacesafterend
}
\newcommand*{\eqrefe}[1]{%
\begingroup
(\color{BrickRed}
\romannumeral-`\x % remove space at the    beginning of \setcitestyle
\setcitestyle{numbers}%
\ref{eq:#1}%
\endgroup
)\ignorespacesafterend
}

\newcommand*{\secrefe}[1]{%
\begingroup
(\color{Aquamarine}
\romannumeral-`\x % remove space at the    beginning of \setcitestyle
\setcitestyle{numbers}%
\ref{#1}%
\endgroup
)\ignorespacesafterend
}

\newtagform{Tags}[\textcolor{BrickRed}]{\color{Black}(}{)}

\begin{document}
\title{Generalized Harmonic Progression}
\date{November 27, 2018}
\author{Jose Risomar Sousa}
\maketitle
\usetagform{Tags}

\begin{abstract}
This paper presents formulae for the sum of the terms of a harmonic progression of order $k$ with integer parameters, $\HP_k(n)$, and for the partial sums of its two associated Fourier series, $C^z_{k}(a,b,n)$ and $S^z_{k}(a,b,n)$. $\HP_k(n)$ is built from the ground up, with a power series for $1/(aj+b)^k$ that is summed over $j$ using Faulhaber's formula. These new formulae are a generalization of the formulae created in a previous paper and were achieved using a slightly modified version of the reasoning employed before.
\end{abstract}

\tableofcontents

\section{Introduction}
Building upon the results of \citena{GHN}, \textit{Generalized Harmonic Numbers}, this new paper demonstrates how to obtain exact formulae for the sum of the first $n$ terms of a harmonic progression of order $k$ with integer parameters, $a$ and $b$: 
\begin{equation} \nonumber
\HP_{k}(n)=\sum_{j=1}^{n} \frac{1}{(a j+b)^k} 
\end{equation}
\indent Even though formulae for $\HP_{k}(n)$ can probably be created using the digamma function, $\psi$, the ones derived here are arguably more interesting -- and also more tractable and easier to work with.\\

We also create formulae for the partial sums of two Fourier series associated with $\HP_{k}(n)$:
\begin{equation} \nonumber
C^z_{k}(a,b,n)=\sum_{j=1}^{n} \frac{1}{(aj+b)^k}\cos{\frac{2\pi(aj+b)}{z}} \text{, and } S^z_{k}(a,b,n)=\sum_{j=1}^{n} \frac{1}{(aj+b)^k}\sin{\frac{2\pi(aj+b)}{z}}
\end{equation}
\indent The first manuscript\citesup{GHN} received some criticism for including too many details, so in this one unnecessary details are omitted, though they can be easily understood from a reading of \citena{GHN}, for the interested reader.\\

This new manuscript is concise and focuses solely on closed-forms for $\HP_{k}(n)$, $C^z_{k}(a,b,n)$ and $S^z_{k}(a,b,n)$\footnote{In a stricter sense, closed-forms do not include integrals.}. Please refer to \citena{Hurwitz} and \citena{Lerch} for the demonstrations on how to obtain the limits of these expressions as $n$ approaches infinity, and to \citena{AC} for how to obtain their analytic continuations.\\

The presented method makes use of Faulhaber's formula\citesup{tres} for the sum of the $i$-th powers of the first $n$ positive integers:
\begin{equation} \nonumber
\sum_{k=1}^{n} {k^i}=\sum_{j=0}^{i} \frac{(-1)^j i! B_j n^{i+1-j}}{(i+1-j)!j!} \text{,} 
\end{equation}
\noindent where $B_j$ are the Bernoulli numbers.\\

Since odd Bernoulli numbers are always zero, except for $B_1$, the above formula can be simplified for even and odd powers as follows:
\begin{equation} \label{eq:soma_pot_par}
\sum_{k=1}^{n} {k^{2i}}=\frac{n^{2i}}{2}+\sum_{j=0}^{i} \frac{(2i)! B_{2j} n^{2i+1-2j}}{(2j)!(2i+1-2j)!} 
\end{equation} 

\begin{equation} \label{eq:soma_pot_impar}
\sum_{k=1}^{n} {k^{2i+1}}=\frac{n^{2i+1}}{2}+\sum_{j=0}^{i} \frac{(2i+1)! B_{2j} n^{2i+2-2j}}{(2j)!(2i+2-2j)!}
\end{equation}

\section{Lagrange's Identities}
In \citena{GHN}, we introduced the indicator function, $k$ divides $n$ ($\mathbbm{1}_{k|n}$), and its analog as key components of the method used to solve the harmonic numbers. For the harmonic progression, it is necessary to modify those functions and obtain their Taylor series -- though unlike in the first paper here we focus on only one of them, as the other one has an analogous behavior.\\

The below is the function that is needed for the odd case. This summation can be written as a closed-form using Lagrange's trigonometric identities:
\begin{equation} \label{eq:k_div_n} \nonumber
\sum_{j=1}^{k}\sin{\frac{2\pi n(a j+b)}{k}}=-\frac{1}{2}\sin{\frac{2\pi b n}{k}}+\frac{1}{2}\sin{2\pi n\left(a+\frac{b}{k}\right)}+\sin{\pi n\left(a+\frac{2b}{k}\right)}\sin{\pi a n}\cot{\frac{\pi a n}{k}} 
\end{equation}
\indent Now, a power series can be derived for the left-hand side of the above equation with the employment of \eqrefe{soma_pot_impar}, and through comparison the following power series can be deduced for each function on the right-hand side (they hold for all real $n$, $a$ and $b$ and for all integer $k\ge 1$):
\begin{multline} \label{eq:k_div_n_series1}
\sum_{i=0}^{\infty}(-1)^i\left(\frac{2\pi b n}{k}\right)^{2i+1}\sum_{j=0}^{i}\left(\frac{(a k/b)^{2j}}{(2j)!(2i+1-2j)!}+\frac{(a k/b)^{2j+1}}{(2j+1)!(2i-2j)!}\right)\\=\sin{2\pi n\left(a+\frac{b}{k}\right)}
\end{multline} 
\begingroup
\small
\begin{multline} \label{eq:k_div_n_series2}
\sum_{i=0}^{\infty}(-1)^i\left(\frac{2\pi b n}{k}\right)^{2i+1}\sum_{j=0}^{i}\left(\frac{(a/b)^{2j}}{(2i+1-2j)!}\sum_{q=0}^{j}\frac{B_{2 q}k^{2j+1-2 q}}{(2j+1-2 q)!(2 q)!}+\frac{(a/b)^{2j+1}}{(2i-2j)!}\sum_{q=0}^{j}\frac{B_{2 q}k^{2j+2-2 q}}{(2j+2-2 q)!(2 q)!}\right)\\=\sin{\pi n\left(a+\frac{2b}{k}\right)}\sin{\pi a n}\cot{\frac{\pi a n}{k}}
\end{multline}
\endgroup
\indent For the harmonic progressions of even order, the analogous function is:
\begin{equation} \nonumber \label{eq:even}
\sum_{j=1}^{k}\cos{\frac{2\pi n(a j+b)}{k}}=-\frac{1}{2}\cos{\frac{2\pi b n}{k}}+\frac{1}{2}\cos{2\pi n\left(a+\frac{b}{k}\right)}+\cos{\pi n\left(a+\frac{2b}{k}\right)}\sin{\pi a n}\cot{\frac{\pi a n}{k}} 
\end{equation}

\section{Approach based on sine} \label{HN_sin_k_pi}
The rationale to build a formula for $\HP_{k}(n)$ is to use the Taylor series expansion of $\sin{2\pi (aj+b)}$,\citesup{Abra} and exploit the fact that it is zero for all integer $aj+b$ (hence the need for $a$ and $b$ to be integers):\\
\begin{equation} \label{eq:start_iden}
\sin{2\pi (a j+b)}=0 \Rightarrow  2\pi(a j+b)=-\sum_{i=1}^{\infty}\frac{(-1)^i}{(2i+1)!}(2\pi (aj+b))^{2i+1} 
\end{equation}\\
\indent If one divides both sides of \eqrefe{start_iden} by $2\pi (a j+b)^2$, a power series for $1/(a j+b)$ that only holds for integer $a j+b$ is obtained.\\

Besides, on the right-hand side of the resulting equation, the exponents of $a j+b$ are positive integers, allowing to apply the Faulhaber formulae from the introduction. This results in a very convoluted power series that can be turned into integrals by means of equations \eqrefe{k_div_n_series1} and \eqrefe{k_div_n_series2}, which in turn were derived using Lagrange's identities. That is a high-level summary of the logic.

\subsection{Harmonic Progression} \label{HN_1_sin_2k_pi}
One starts by dividing both sides of \eqrefe{start_iden} by $2\pi (a j+b)^2$:
\begin{equation} \label{eq:start_iden_2}
\frac{1}{a j+b}=\sum_{i=0}^{\infty}\frac{(-1)^i (2\pi)^{2i+2}(a j+b)^{2i+1}}{(2i+3)!} 
\end{equation}\\
\indent Now, $1/(a j+b)$ is summed over $j$, $(a j+b)^{2i+1}$ is expanded with the binomial theorem and the Faulhaber formulae, from \eqrefe{soma_pot_par} and \eqrefe{soma_pot_impar}, are applied, yielding the following power series after all the calculations are done:
\begingroup
\small
\begin{multline} \nonumber
\sum_{j=1}^{n}\frac{1}{a j+b}=-\frac{1}{2b}+\frac{1}{2}\sum_{i=0}^{\infty}(-1)^i\frac{(2\pi b)^{2i+2}(2i+1)!}{b(2i+3)!}\sum_{j=0}^{i}\left(\frac{(a n/b)^{2j}}{(2j)!(2i+1-2j)!}+\frac{(a n/b)^{2j+1}}{(2j+1)!(2i-2j)!}\right)+\\ 
+\sum_{i=0}^{\infty}(-1)^i\frac{(2\pi b)^{2i+2}(2i+1)!}{b(2i+3)!}\sum_{j=0}^{i}\left(\frac{(a/b)^{2j}}{(2i+1-2j)!}\sum_{q=0}^{j}\frac{B_{2 q}n^{2j+1-2 q}}{(2j+1-2 q)!(2 q)!}+\frac{(a/b)^{2j+1}}{(2i-2j)!}\sum_{q=0}^{j}\frac{B_{2 q}n^{2j+2-2 q}}{(2j+2-2 q)!(2 q)!}\right) 
\end{multline}
\endgroup\\
\indent The above sums can be manipulated conveniently and then obtained from \eqrefe{k_div_n_series1} and \eqrefe{k_div_n_series2}, the two equations derived previously, by means of decompositions into linear combinations followed by integrations. Though the details were omitted here, the reader can refer to the precursor paper\citesup{GHN} for a more detailed description of the steps involved.\\

After all the appropriate calculations are performed, one ends up with: 
\begingroup
\small
\begin{equation} \nonumber
\frac{1}{2b}\sum_{i=0}^{\infty}(-1)^i\frac{(2\pi b)^{2i+2}(2i+1)!}{(2i+3)!}\sum_{j=0}^{i}\left(\frac{(a n/b)^{2j}}{(2j)!(2i+1-2j)!}+\frac{(a n/b)^{2j+1}}{(2j+1)!(2i-2j)!}\right)=\frac{2\pi(a n+b)-\sin{2\pi(a n+b)}}{4\pi(a n+b)^2} 
\end{equation} 
\begin{multline} \nonumber
\sum_{i=0}^{\infty}(-1)^i\frac{(2\pi b)^{2i+2}(2i+1)!}{b(2i+3)!}\sum_{j=0}^{i}\left(\frac{(a/b)^{2j}}{(2i+1-2j)!}\sum_{q=0}^{j}\frac{B_{2 q}n^{2j+1-2 q}}{(2j+1-2 q)!(2 q)!}+\frac{(a/b)^{2j+1}}{(2i-2j)!}\sum_{q=0}^{j}\frac{B_{2 q}n^{2j+2-2 q}}{(2j+2-2 q)!(2 q)!}\right)\\
=2\pi\int_{0}^{1}(1-u)\sin{\pi a n u}\sin{\pi(a n+2b)u}\cot{\pi a u}\,du 
\end{multline}
\endgroup\\
\indent Now, adding up all the results (disregarding the sine of multiples of $\pi$), one arrives at a formula for $\HP(n)$:
\begin{equation} \nonumber
\sum_{j=1}^{n}\frac{1}{a j+b}=-\frac{1}{2b}+\frac{1}{2(a n+b)}+ 2\pi\int_{0}^{1}(1-u)\sin{\pi(a n+2b)u}\sin{\pi a n u}\cot{\pi a u}\,du 
\end{equation}
\indent In the next sections a brief generalization of this result is stated.

\subsection{General Formula} \label{Gen_sin_2k_pi}

If one keeps dividing \eqrefe{start_iden_2} by $a j+b$, similar recursions to the ones obtained for the generalized harmonic numbers are obtained. All the results presented next follow from reasonings analogous to those from reference \citena{GHN}.

\subsection{Harmonic Progression of Order $2k$} \label{HN_2k_sin_2k_pi}
One has the following recursion for $\HP_{2k}(n)$:
\begin{multline} \nonumber
\HP_{2k}(n)=-\frac{1}{2 b^{2k}}\sum_{j=0}^{k}\frac{(-1)^{j}(2\pi b)^{2j}}{(2j+1)!}+\frac{1}{2 (an+b)^{2k}}\sum_{j=0}^{k}\frac{(-1)^{j}(2\pi(an+b))^{2j}}{(2j+1)!}
\\-\sum_{j=1}^{k-1}\frac{(-1)^{k-j}(2\pi)^{2k-2j}}{(2k+1-2j)!}\HP_{2j}(n)-\frac{(-1)^k(2\pi)^{2k}}{(2k)!}\int_{0}^{1}(1-u)^{2k}\cos{\pi(a n+2b)u}\sin{\pi a n u}\cot{\pi a u}\,du  \text{}
\end{multline}
\indent Note that $\HP_{0}(n)=0$ for all positive integer $n$ (just like $H_0(n)=0$, previously). Therefore, from the recursion, it follows that for all integer $a\neq 0$, $b\neq 0$ and $k\geq 1$:
\begin{multline} \nonumber
\sum_{j=1}^{n}\frac{1}{(a j+b)^{2k}}=-\frac{1}{2b^{2k}}+\frac{1}{2(a n+b)^{2k}}+\\-(-1)^{k}(2\pi)^{2k}\int_{0}^{1}\sum_{j=0}^{k}\frac{B_{2j}\left(2-2^{2j}\right)(1-u)^{2k-2j}}{(2j)!(2k-2j)!}\cos{\pi(a n+2b)u}\sin{\pi a n u}\cot{\pi a u}\,du 
\end{multline}
\indent And if the product of cosine and sine is turned into a sum of sines, and if $1-u$ is replaced with $u$ (this does not change the integral or each of its individual parts), a different, perhaps more useful, way to express this formula is obtained:
\begin{multline} \label{eq:HP_2k}
\sum_{j=1}^{n}\frac{1}{(a j+b)^{2k}}=-\frac{1}{2b^{2k}}+\frac{1}{2(a n+b)^{2k}}\\-\frac{(-1)^{k}(2\pi)^{2k}}{2}\int_{0}^{1}\sum_{j=0}^{k}\frac{B_{2j}\left(2-2^{2j}\right)u^{2k-2j}}{(2j)!(2k-2j)!}\left(\sin{2\pi(a n+b)u}-\sin{2\pi b u}\right)\cot{\pi a u}\,du 
\end{multline}
\indent This formula also holds for the generalized harmonic numbers ($a=1,b=0$), if the term $-1/(2b^{2k})$ is disregarded. In fact, if any term of the equation that has a null denominator is disregarded, the equation still holds.

\subsection{Harmonic Progression of Order $2k+1$} \label{HN_2k+1_sin_2k_pi}

One has the following recursion for $\HP_{2k+1}(n)$:
\begin{multline} \nonumber
\HP_{2k+1}(n)=-\frac{1}{2 b^{2k+1}}\sum_{j=0}^{k}\frac{(-1)^{j}(2\pi b)^{2j}}{(2j+1)!}+\frac{1}{2 (an+b)^{2k+1}}\sum_{j=0}^{k}\frac{(-1)^{j}(2\pi(an+b))^{2j}}{(2j+1)!}
\\-\sum_{j=0}^{k-1}\frac{(-1)^{k-j}(2\pi)^{2k-2j}}{(2k+1-2j)!}\HP_{2j+1}(n)+\frac{(-1)^k(2\pi)^{2k+1}}{(2k+1)!}\int_{0}^{1}(1-u)^{2k+1}\sin{\pi(a n+2b)u}\sin{\pi a n u}\cot{\pi a u}\,du \text{} 
\end{multline}
\indent Therefore, for all integer $a\neq 0$, $b\neq 0$ and $k\geq 0$:
\begin{multline} \nonumber
\sum_{j=1}^{n}\frac{1}{(a j+b)^{2k+1}}=-\frac{1}{2b^{2k+1}}+\frac{1}{2(a n+b)^{2k+1}}+\\
+(-1)^{k}(2\pi)^{2k+1}\int_{0}^{1}\sum_{j=0}^{k}\frac{B_{2j}\left(2-2^{2j}\right)(1-u)^{2k+1-2j}}{(2j)!(2k+1-2j)!}\sin{\pi(a n+2b)u}\sin{\pi a nu}\cot{\pi au}\,du 
\end{multline}
\indent Again, if the product of sines is turned into a sum of cosines, and if $1-u$ is replaced with $u$ (this only flips the sign of the integral, ditto for each of its individual parts), the below is obtained:
\begin{multline} \label{eq:HP_2k+1}
\sum_{j=1}^{n}\frac{1}{(a j+b)^{2k+1}}=-\frac{1}{2b^{2k+1}}+\frac{1}{2(a n+b)^{2k+1}}\\
+\frac{(-1)^{k}(2\pi)^{2k+1}}{2}\int_{0}^{1}\sum_{j=0}^{k}\frac{B_{2j}\left(2-2^{2j}\right)u^{2k+1-2j}}{(2j)!(2k+1-2j)!}\left(\cos{2\pi(a n+b)u}-\cos{2\pi bu}\right)\cot{\pi au}\,du 
\end{multline}

\subsection{Generating functions}
As seen in \citena{GHN}, the polynomials in $u$ within the integrals are generated by the functions:
\begin{equation} \nonumber
f(x)=\frac{x\cos{x(1-u)}}{\sin{x}} \Rightarrow \frac{f^{(2k)}(0)}{(2k)!}=(-1)^k\sum_{j=0}^{k}\frac{B_{2j}\left(2-2^{2j}\right)(1-u)^{2k-2j}}{(2j)!(2k-2j)!} \text{, and}
\end{equation}
\begin{equation} \nonumber
f(x)=\frac{x\sin{x(1-u)}}{\sin{x}} \Rightarrow \frac{f^{(2k+1)}(0)}{(2k+1)!}=(-1)^k\sum_{j=0}^{k}\frac{B_{2j}\left(2-2^{2j}\right)(1-u)^{2k+1-2j}}{(2j)!(2k+1-2j)!}
\end{equation}

\section{Approach based on exponential}
This approach consists in using equation $e^{2\pi\ii(aj+b)}=1$, which combines the cosine and sine equations into one. It has a possible advantage over the sine-based approach, since it yields a single formula for both odd and even powers. For this one, we skip the step-by-step demonstration and go straight to the final result.

\subsection{General formula}
For all integer $a$, $b$ and $k\ge 1$:
\begin{multline} \nonumber
\sum_{j=1}^{n}\frac{1}{(a j+b)^{k}}=-\frac{1}{2b^{k}}+\frac{1}{2(a n+b)^{k}}\\+\frac{\ii(2\pi\ii)^{k}}{2}\int_{0}^{1}\sum_{j=0}^{k}\frac{B_{j}(1-u)^{k-j}}{j!(k-j)!}\left(e^{2\pi\ii(a n+b)u}-e^{2\pi\ii b\,u}\right)\cot{\pi a u}\,du 
\end{multline}
\indent A constant $c$ can be introduced into the formula such that, provided that $c\,a$ and $c\,b$ are both integers, the following modified formula still holds:
\begin{multline} \nonumber
\sum_{j=1}^{n}\frac{1}{(a j+b)^{k}}=-\frac{1}{2b^{k}}+\frac{1}{2(a n+b)^{k}}\\+\frac{\ii(-2\pi\ii\,c)^{k}}{2}\int_{0}^{1}\sum_{j=0}^{k}\frac{B_{j}u^{k-j}}{j!(k-j)!}\left(e^{2\pi\ii\,c(a n+b)u}-e^{2\pi\ii\,c\,b\,u}\right)\cot{\pi c\,a u}\,du 
\end{multline}
\indent Note the negative sign stems from a change of variables (replacing $1-u$ with $u$), it does not come from the introduced constant. The constant can be useful. For example, if $a=1$ and $b=1/3$, setting $c=3$ makes the formula right.

\subsection{Decoupling}
It is possible to remove the complex numbers out of the picture, in which case the formula becomes:
\begin{multline} \nonumber
\sum_{j=1}^{n}\frac{1}{(a j+b)^{k}}=-\frac{1}{2b^{k}}+\frac{1}{2(a n+b)^{k}}\\-(2\pi)^{k}\int_{0}^{1}\sum_{j=0}^{k}\frac{B_{j}(1-u)^{k-j}}{j!(k-j)!}\cos{\left((a n+2b)\pi u+\frac{k\pi}{2}\right)}\sin{\pi a n u}\cot{\pi a u}\,du \text{,}
\end{multline}
\noindent which in turn can be transformed into forms quite similar to those from section \secrefe{Gen_sin_2k_pi}.\\

For the even powers one has:
\begin{multline} \nonumber %\label{eq:dec_HP_2k}
\sum_{j=1}^{n}\frac{1}{(a j+b)^{2k}}=-\frac{1}{2b^{2k}}+\frac{1}{2(a n+b)^{2k}}\\-\frac{(-1)^{k}(2\pi)^{2k}}{2}\int_{0}^{1}\left(-\frac{u^{2k-1}}{2(2k-1)!}+\sum_{j=0}^{k}\frac{B_{2j}u^{2k-2j}}{(2j)!(2k-2j)!}\right)\left(\sin{2\pi(a n+b)u}-\sin{2\pi b u}\right)\cot{\pi a u}\,du \text{,}
\end{multline}
\noindent and for the odd powers:
\begin{multline} \nonumber
\sum_{j=1}^{n}\frac{1}{(a j+b)^{2k+1}}=-\frac{1}{2b^{2k+1}}+\frac{1}{2(a n+b)^{2k+1}}\\
+\frac{(-1)^{k}(2\pi)^{2k+1}}{2}\int_{0}^{1}\left(-\frac{u^{2k}}{2(2k)!}+\sum_{j=0}^{k}\frac{B_{2j}u^{2k+1-2j}}{(2j)!(2k+1-2j)!}\right)\left(\cos{2\pi(a n+b)u}-\cos{2\pi bu}\right)\cot{\pi au}\,du 
\end{multline}

\subsection{Generating function}
The generating function of the polynomial in $u$ within the integral is:
\begin{equation} \nonumber
f(x)=\frac{x\, e^{x(1-u)}}{e^x-1} \Rightarrow \frac{f^{(k)}(0)}{k!}=\sum_{j=0}^{k}\frac{B_{j}(1-u)^{k-j}}{j!(k-j)!} 
\end{equation}

\section{The partial Fourier series} \label{Final}
Let us briefly recall the formulae that were found for $C^z_{k}(n)$ and $S^z_{k}(n)$, the partial sums of the Fourier series associated with the generalized harmonic numbers, $H_k(n)$, from reference \citena{GHN}. Next to each one, their harmonic progression analogs, $C^z_{k}(a,b,n)$ and $S^z_{k}(a,b,n)$, are displayed, which are based entirely on analogy.\\

The following expressions hold for all complex $z$, $a$ and $b$ (unlike the $\HP_k(n)$ formulae), and for all integer $n\geq 1$. As mentioned before, any term that has zero in the denominator can be disregarded, and the equation still holds (technically, the limit as the parameter tends to zero is taken). See \citena{Lerch} and \citena{AC}, for examples.\\ 

By definition $H_0(n)=0$ and $\HP_{0}(n)=0$ for all positive integer $n$, so they have no effect in the summations and, to avoid confusion, they are skipped in the formulae.

\subsection{$C^z_{2k}(n)$ and $C^z_{2k}(a,b,n)$} \label{Final_1}
For all integer $k \geq 1$:
\begingroup
\footnotesize
\begin{multline} \nonumber
\sum_{j=1}^{n}\frac{1}{j^{2k}}\cos{\frac{2\pi j}{z}}=\frac{1}{2n^{2k}}\left(\cos{\frac{2\pi n}{z}}-\sum_{j=0}^{k}\frac{(-1)^j (\frac{2\pi n}{z})^{2j}}{(2j)!}\right)+\sum_{j=1}^{k}\frac{(-1)^{k-j}(\frac{2\pi}{z})^{2k-2j}}{(2k-2j)!}H_{2j}(n)
\\+\frac{(-1)^k(\frac{2\pi}{z})^{2k}}{2(2k-1)!}\int_{0}^{1}(1-u)^{2k-1}\sin{\frac{2\pi n u}{z}}\cot{\frac{\pi u}{z}}\,du 
\end{multline}

\begin{multline} \nonumber
\sum_{j=1}^{n}\frac{1}{(a j+b)^{2k}}\cos{\frac{2\pi(a j+b)}{z}}=-\frac{1}{2b^{2k}}\left(\cos{\frac{2\pi b}{z}}-\sum_{j=0}^{k}\frac{(-1)^j (\frac{2\pi b}{z})^{2j}}{(2j)!}\right)
\\+\frac{1}{2(a n+b)^{2k}}\left(\cos{\frac{2\pi(a n+b)}{z}}-\sum_{j=0}^{k}\frac{(-1)^j (\frac{2\pi(a n+b)}{z})^{2j}}{(2j)!}\right)+\sum_{j=1}^{k}\frac{(-1)^{k-j}(\frac{2\pi}{z})^{2k-2j}}{(2k-2j)!}\HP_{2j}(n)
\\+\frac{(-1)^k(\frac{2\pi}{z})^{2k}}{2(2k-1)!}\int_{0}^{1}(1-u)^{2k-1}\left(\sin{\frac{2\pi(a n+b)u}{z}}-\sin{\frac{2\pi b u}{z}}\right)\cot{\frac{\pi a u}{z}}\,du 
\end{multline}
\endgroup 

\subsection{$S^z_{2k+1}(n)$ and $S^z_{2k+1}(a,b,n)$} \label{Final_2}
For all integer $k \geq 0$:
\begingroup
\footnotesize
\begin{multline} \nonumber
\sum_{j=1}^{n}\frac{1}{j^{2k+1}}\sin{\frac{2\pi j}{z}}=\frac{1}{2n^{2k+1}}\left(\sin{\frac{2\pi n}{z}}-\sum_{j=0}^{k}\frac{(-1)^j (\frac{2\pi n}{z})^{2j+1}}{(2j+1)!}\right)+\sum_{j=1}^{k}\frac{(-1)^{k-j}(\frac{2\pi}{z})^{2k+1-2j}}{(2k+1-2j)!}H_{2j}(n)\\+\frac{(-1)^k(\frac{2\pi}{z})^{2k+1}}{2(2k)!}\int_{0}^{1}(1-u)^{2k}\sin{\frac{2\pi n u}{z}}\cot{\frac{\pi u}{z}}\,du 
\end{multline}

\begin{multline} \nonumber
\sum_{j=1}^{n}\frac{1}{(a j+b)^{2k+1}}\sin{\frac{2\pi(a j+b)}{z}}=-\frac{1}{2b^{2k+1}}\left(\sin{\frac{2\pi b}{z}}-\sum_{j=0}^{k}\frac{(-1)^j (\frac{2\pi b}{z})^{2j+1}}{(2j+1)!}\right)
\\+\frac{1}{2(a n+b)^{2k+1}}\left(\sin{\frac{2\pi(a n+b)}{z}}-\sum_{j=0}^{k}\frac{(-1)^j (\frac{2\pi(a n+b)}{z})^{2j+1}}{(2j+1)!}\right)+\sum_{j=1}^{k}\frac{(-1)^{k-j}(\frac{2\pi}{z})^{2k+1-2j}}{(2k+1-2j)!}\HP_{2j}(n)
\\+\frac{(-1)^k(\frac{2\pi}{z})^{2k+1}}{2(2k)!}\int_{0}^{1}(1-u)^{2k}\left(\sin{\frac{2\pi(a n+b)u}{z}}-\sin{\frac{2\pi b u}{z}}\right)\cot{\frac{\pi a u}{z}}\,du 
\end{multline}
\endgroup 

\subsection{$C^z_{2k+1}(n)$ and $C^z_{2k+1}(a,b,n)$} \label{Final_3}
For all integer $k \geq 0$:
\begingroup
\footnotesize
\begin{multline} \nonumber
\sum_{j=1}^{n}\frac{1}{j^{2k+1}}\cos{\frac{2\pi j}{z}}=\frac{1}{2n^{2k+1}}\left(\cos{\frac{2\pi n}{z}}-\sum_{j=0}^{k}\frac{(-1)^j (\frac{2\pi n}{z})^{2j}}{(2j)!}\right)+\sum_{j=0}^{k}\frac{(-1)^{k-j}(\frac{2\pi}{z})^{2k-2j}}{(2k-2j)!}H_{2j+1}(n)\\+\frac{(-1)^k(\frac{2\pi}{z})^{2k+1}}{2(2k)!}\int_{0}^{1}(1-u)^{2k}\left(\cos{\frac{2\pi n u}{z}}-1\right)\cot{\frac{\pi u}{z}}\,du 
\end{multline}

\begin{multline} \nonumber
\sum_{j=1}^{n}\frac{1}{(a j+b)^{2k+1}}\cos{\frac{2\pi(a j+b)}{z}}=-\frac{1}{2b^{2k+1}}\left(\cos{\frac{2\pi b}{z}}-\sum_{j=0}^{k}\frac{(-1)^j (\frac{2\pi b}{z})^{2j}}{(2j)!}\right)
\\+\frac{1}{2(a n+b)^{2k+1}}\left(\cos{\frac{2\pi(a n+b)}{z}}-\sum_{j=0}^{k}\frac{(-1)^j (\frac{2\pi(a n+b)}{z})^{2j}}{(2j)!}\right)+\sum_{j=0}^{k}\frac{(-1)^{k-j}(\frac{2\pi}{z})^{2k-2j}}{(2k-2j)!}\HP_{2j+1}(n)
\\+\frac{(-1)^k(\frac{2\pi}{z})^{2k+1}}{2(2k)!}\int_{0}^{1}(1-u)^{2k}\left(\cos{\frac{2\pi(a n+b)u}{z}}-\cos{\frac{2\pi b u}{z}}\right)\cot{\frac{\pi a u}{z}}\,du 
\end{multline}
\endgroup 

\subsection{$S^z_{2k}(n)$ and $S^z_{2k}(a,b,n)$} \label{Final_4}
For all integer $k \geq 1$:
\begingroup
\footnotesize
\begin{multline} \nonumber
\sum_{j=1}^{n}\frac{1}{j^{2k}}\sin{\frac{2\pi j}{z}}=\frac{1}{2n^{2k}}\left(\sin{\frac{2\pi n}{z}}-\sum_{j=0}^{k-1}\frac{(-1)^j (\frac{2\pi n}{z})^{2j+1}}{(2j+1)!}\right)-\sum_{j=0}^{k-1}\frac{(-1)^{k-j}(\frac{2\pi}{z})^{2k-1-2j}}{(2k-1-2j)!}H_{2j+1}(n)\\
-\frac{(-1)^k(\frac{2\pi}{z})^{2k}}{2(2k-1)!}\int_{0}^{1}(1-u)^{2k-1}\left(\cos{\frac{2\pi n u}{z}}-1\right)\cot{\frac{\pi u}{z}}\,du 
\end{multline}

\begin{multline} \nonumber
\sum_{j=1}^{n}\frac{1}{(a j+b)^{2k}}\sin{\frac{2\pi(a j+b)}{z}}=-\frac{1}{2b^{2k}}\left(\sin{\frac{2\pi b}{z}}-\sum_{j=0}^{k-1}\frac{(-1)^j (\frac{2\pi b}{z})^{2j+1}}{(2j+1)!}\right)
\\+\frac{1}{2(a n+b)^{2k}}\left(\sin{\frac{2\pi(a n+b)}{z}}-\sum_{j=0}^{k-1}\frac{(-1)^j (\frac{2\pi(a n+b)}{z})^{2j+1}}{(2j+1)!}\right)-\sum_{j=0}^{k-1}\frac{(-1)^{k-j}(\frac{2\pi}{z})^{2k-1-2j}}{(2k-1-2j)!}\HP_{2j+1}(n)
\\-\frac{(-1)^k(\frac{2\pi}{z})^{2k}}{2(2k-1)!}\int_{0}^{1}(1-u)^{2k-1}\left(\cos{\frac{2\pi(a n+b)u}{z}}-\cos{\frac{2\pi b u}{z}}\right)\cot{\frac{\pi a u}{z}}\,du 
\end{multline} 
\endgroup

\end{document}